\def\marker{\>\hbox{${\vcenter{\vbox{
    \hrule height 0.4pt\hbox{\vrule width 0.4pt height 6pt
    \kern6pt\vrule width 0.4pt}\hrule height 0.4pt}}}$}\>}
\def\gpic#1{#1
     \medskip\par\noindent{\centerline{\box\graph}} \medskip}
\documentclass{article}
\usepackage{amsfonts}
\usepackage{pslatex}

\author{Daniel Cranston
\footnote{University of Illinois, Urbana-Champaign.  This work was supported by the Mathematical, Information, and 
Computational Sciences Division subprogram of the Office of Advanced 
Scientific Computing Research, Office of Science, U.S. Department of 
Energy, under Contract W-31-109-ENG-38.} \\
\texttt{dcransto@uiuc.edu}
}

\title{A Strong Edge-Coloring of Graphs with Maximum Degree 4 Using 22 Colors
}

\newtheorem{lemma}{Lemma}
\newtheorem{theorem}[lemma]{Theorem}

\begin{document}
\pagestyle{headings}

%\begin{titlepage}
\date{May 3, 2005}
\maketitle
\begin{abstract}
In 1985, Erd\H{o}s and Ne\'{s}etril conjectured that the strong edge-coloring number of a graph is bounded above by $\frac{5}{4}\Delta^2$ when $\Delta$ is even and $\frac{1}{4}(5\Delta^2-2\Delta+1)$ when $\Delta$ is odd.  They gave a simple construction which requires this many colors.  The conjecture has been verified for $\Delta\leq 3$.  For $\Delta=4$, the conjectured bound is 20. Previously, the best known upper bound was 23 due to Horak.  In this paper we give an algorithm that uses at most 22 colors.
\end{abstract}
%\end{titlepage}

\section{Introduction}

A {\it proper edge-coloring} is an assignment of a color to each edge of a graph so that no two edges with a common endpoint receive the same color.  A {\it strong edge-coloring} is a proper edge-coloring, with the further condition that no two edges with the same color lie on a path of length three.  The {\it strong edge chromatic number} is the minimum number of colors that allow a strong edge-coloring.  In this paper we consider the maximum possible strong edge chromatic number as a function of the maximum degree of the graph.  For other variations of the problem, we refer the reader to a brief survey by West \cite{We} and a paper by Faudree, Schelp, Gy\'{a}rf\'{a}s and Tuza \cite{FGST}.

We use $\Delta$ to denote the maximum degree of the graph.  In 1985 Erd\H{o}s and Ne\v{s}et\v{r}il conjectured that the strong edge chromatic number of a graph is at most $\frac{5}{4}\Delta^2$ for $\Delta$ even and $\frac{1}{4}(5\Delta^2-2\Delta+1)$ for $\Delta$ odd; they gave a construction that showed this number is necessary.  Andersen proved the conjecture for the case $\Delta = 3$ \cite{An}.  In this paper, we consider the case $\Delta=4$.  

Erd\H{o}s and Ne\v{s}et\v{r}il's construction for $\Delta=4$ is shown in figure (1a).  To form this graph, begin with a 5-cycle, then expand each vertex into two nonadjacent vertices who inherit all the neighbors of the original vertex.  The graph has 20 edges, and contains no induced $2K_2$ (in fact, this is the largest graph that contains no induced $2K_2$ and has $\Delta=4$ \cite{Ch}).  Hence, in a strong edge-coloring, every edge must receive its own color.  The best upper bound previously known was 23 colors, proven by Horak \cite{Ho}; we improve this upper bound to 22 colors.

We refer to the color classes as the integers from 1 to 22.  A greedy coloring algorithm sequentially colors the edges, using the least color class that is not already prohibited from use on an edge at the time the edge is colored.  By the {\it neighborhood} of an edge, we mean the edges which are distance at most 1 from the edge.  Intuitively, this is the set of edges whose color could potentially restrict the color of that edge.   We use the notation $N(e)$ to mean the edges in the neighborhood of $e$ that are colored before edge $e$.  Figure (1b) shows that the neighborhood of an edge has size at most 24.  So for every edge order, we have $|N(e)|\leq 24$ for each edge.  Thus, for every edge order, the greedy algorithm produces a strong edge-coloring that uses at most 25 colors.
However, there is always some order of the edges for which the greedy algorithm uses exactly the minimum number of colors required.  Our aim in this paper is to construct an order of the edges such that the greedy algorithm uses at most 22 colors.  Throughout this paper, when we use the term coloring, we mean strong edge-coloring.  Each connected component of $G$ can be colored independently of other components, so we assume $G$ is connected.  We allow our graphs to include loops and multiple edges.  We use $\delta$ to denote the minimum degree of the graph and we use $d(v)$ to denote the degree of vertex $v$.  The girth of a graph is the length of the shortest cycle.

\begin{figure}[hbtp]
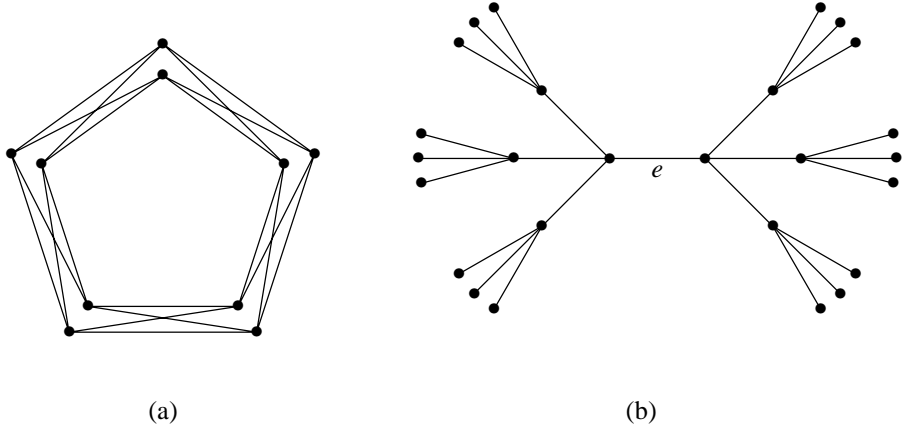

  \vspace{4pt}

\gpic{
\expandafter\ifx\csname graph\endcsname\relax \csname newbox\endcsname\graph\fi
\expandafter\ifx\csname graphtemp\endcsname\relax \csname newdimen\endcsname\graphtemp\fi
\setbox\graph=\vtop{\vskip 0pt\hbox{%
    \graphtemp=.5ex\advance\graphtemp by 0.270in
    \rlap{\kern 0.877in\lower\graphtemp\hbox to 0pt{\hss $\bullet$\hss}}%
    \graphtemp=.5ex\advance\graphtemp by 0.847in
    \rlap{\kern 1.671in\lower\graphtemp\hbox to 0pt{\hss $\bullet$\hss}}%
    \graphtemp=.5ex\advance\graphtemp by 1.780in
    \rlap{\kern 1.368in\lower\graphtemp\hbox to 0pt{\hss $\bullet$\hss}}%
    \graphtemp=.5ex\advance\graphtemp by 1.780in
    \rlap{\kern 0.386in\lower\graphtemp\hbox to 0pt{\hss $\bullet$\hss}}%
    \graphtemp=.5ex\advance\graphtemp by 0.847in
    \rlap{\kern 0.083in\lower\graphtemp\hbox to 0pt{\hss $\bullet$\hss}}%
    \graphtemp=.5ex\advance\graphtemp by 0.437in
    \rlap{\kern 0.877in\lower\graphtemp\hbox to 0pt{\hss $\bullet$\hss}}%
    \graphtemp=.5ex\advance\graphtemp by 0.899in
    \rlap{\kern 1.512in\lower\graphtemp\hbox to 0pt{\hss $\bullet$\hss}}%
    \graphtemp=.5ex\advance\graphtemp by 1.645in
    \rlap{\kern 1.270in\lower\graphtemp\hbox to 0pt{\hss $\bullet$\hss}}%
    \graphtemp=.5ex\advance\graphtemp by 1.645in
    \rlap{\kern 0.485in\lower\graphtemp\hbox to 0pt{\hss $\bullet$\hss}}%
    \graphtemp=.5ex\advance\graphtemp by 0.899in
    \rlap{\kern 0.242in\lower\graphtemp\hbox to 0pt{\hss $\bullet$\hss}}%
    \special{pn 8}%
    \special{pa 877 437}%
    \special{pa 1671 847}%
    \special{pa 877 270}%
    \special{pa 1512 899}%
    \special{pa 877 437}%
    \special{fp}%
    \special{pa 1270 1645}%
    \special{pa 1671 847}%
    \special{pa 1368 1780}%
    \special{pa 1512 899}%
    \special{pa 1270 1645}%
    \special{fp}%
    \special{pa 1270 1645}%
    \special{pa 386 1780}%
    \special{pa 1368 1780}%
    \special{pa 485 1645}%
    \special{pa 1270 1645}%
    \special{fp}%
    \special{pa 242 899}%
    \special{pa 386 1780}%
    \special{pa 83 847}%
    \special{pa 485 1645}%
    \special{pa 242 899}%
    \special{fp}%
    \special{pa 877 437}%
    \special{pa 83 847}%
    \special{pa 877 270}%
    \special{pa 242 899}%
    \special{pa 877 437}%
    \special{fp}%
    \graphtemp=.5ex\advance\graphtemp by 2.207in
    \rlap{\kern 0.877in\lower\graphtemp\hbox to 0pt{\hss (a)\hss}}%
    \special{pa 3214 871}%
    \special{pa 3715 871}%
    \special{fp}%
    \graphtemp=.5ex\advance\graphtemp by 0.938in
    \rlap{\kern 3.465in\lower\graphtemp\hbox to 0pt{\hss $e$\hss}}%
    \graphtemp=.5ex\advance\graphtemp by 0.871in
    \rlap{\kern 3.214in\lower\graphtemp\hbox to 0pt{\hss $\bullet$\hss}}%
    \graphtemp=.5ex\advance\graphtemp by 0.871in
    \rlap{\kern 3.715in\lower\graphtemp\hbox to 0pt{\hss $\bullet$\hss}}%
    \graphtemp=.5ex\advance\graphtemp by 0.517in
    \rlap{\kern 4.069in\lower\graphtemp\hbox to 0pt{\hss $\bullet$\hss}}%
    \graphtemp=.5ex\advance\graphtemp by 0.871in
    \rlap{\kern 4.216in\lower\graphtemp\hbox to 0pt{\hss $\bullet$\hss}}%
    \graphtemp=.5ex\advance\graphtemp by 1.225in
    \rlap{\kern 4.069in\lower\graphtemp\hbox to 0pt{\hss $\bullet$\hss}}%
    \special{pa 4069 517}%
    \special{pa 3715 871}%
    \special{fp}%
    \special{pa 4216 871}%
    \special{pa 3715 871}%
    \special{fp}%
    \special{pa 4069 1225}%
    \special{pa 3715 871}%
    \special{fp}%
    \graphtemp=.5ex\advance\graphtemp by 1.225in
    \rlap{\kern 2.860in\lower\graphtemp\hbox to 0pt{\hss $\bullet$\hss}}%
    \graphtemp=.5ex\advance\graphtemp by 0.871in
    \rlap{\kern 2.713in\lower\graphtemp\hbox to 0pt{\hss $\bullet$\hss}}%
    \graphtemp=.5ex\advance\graphtemp by 0.517in
    \rlap{\kern 2.860in\lower\graphtemp\hbox to 0pt{\hss $\bullet$\hss}}%
    \special{pa 2860 1225}%
    \special{pa 3214 871}%
    \special{fp}%
    \special{pa 2713 871}%
    \special{pa 3214 871}%
    \special{fp}%
    \special{pa 2860 517}%
    \special{pa 3214 871}%
    \special{fp}%
    \graphtemp=.5ex\advance\graphtemp by 0.083in
    \rlap{\kern 4.319in\lower\graphtemp\hbox to 0pt{\hss $\bullet$\hss}}%
    \graphtemp=.5ex\advance\graphtemp by 0.163in
    \rlap{\kern 4.423in\lower\graphtemp\hbox to 0pt{\hss $\bullet$\hss}}%
    \graphtemp=.5ex\advance\graphtemp by 0.267in
    \rlap{\kern 4.503in\lower\graphtemp\hbox to 0pt{\hss $\bullet$\hss}}%
    \special{pa 4319 83}%
    \special{pa 4069 517}%
    \special{fp}%
    \special{pa 4423 163}%
    \special{pa 4069 517}%
    \special{fp}%
    \special{pa 4503 267}%
    \special{pa 4069 517}%
    \special{fp}%
    \graphtemp=.5ex\advance\graphtemp by 0.742in
    \rlap{\kern 4.699in\lower\graphtemp\hbox to 0pt{\hss $\bullet$\hss}}%
    \graphtemp=.5ex\advance\graphtemp by 0.871in
    \rlap{\kern 4.717in\lower\graphtemp\hbox to 0pt{\hss $\bullet$\hss}}%
    \graphtemp=.5ex\advance\graphtemp by 1.001in
    \rlap{\kern 4.699in\lower\graphtemp\hbox to 0pt{\hss $\bullet$\hss}}%
    \special{pa 4699 742}%
    \special{pa 4216 871}%
    \special{fp}%
    \special{pa 4717 871}%
    \special{pa 4216 871}%
    \special{fp}%
    \special{pa 4699 1001}%
    \special{pa 4216 871}%
    \special{fp}%
    \graphtemp=.5ex\advance\graphtemp by 1.476in
    \rlap{\kern 4.503in\lower\graphtemp\hbox to 0pt{\hss $\bullet$\hss}}%
    \graphtemp=.5ex\advance\graphtemp by 1.579in
    \rlap{\kern 4.423in\lower\graphtemp\hbox to 0pt{\hss $\bullet$\hss}}%
    \graphtemp=.5ex\advance\graphtemp by 1.659in
    \rlap{\kern 4.319in\lower\graphtemp\hbox to 0pt{\hss $\bullet$\hss}}%
    \special{pa 4503 1476}%
    \special{pa 4069 1225}%
    \special{fp}%
    \special{pa 4423 1579}%
    \special{pa 4069 1225}%
    \special{fp}%
    \special{pa 4319 1659}%
    \special{pa 4069 1225}%
    \special{fp}%
    \graphtemp=.5ex\advance\graphtemp by 1.659in
    \rlap{\kern 2.610in\lower\graphtemp\hbox to 0pt{\hss $\bullet$\hss}}%
    \graphtemp=.5ex\advance\graphtemp by 1.579in
    \rlap{\kern 2.506in\lower\graphtemp\hbox to 0pt{\hss $\bullet$\hss}}%
    \graphtemp=.5ex\advance\graphtemp by 1.476in
    \rlap{\kern 2.426in\lower\graphtemp\hbox to 0pt{\hss $\bullet$\hss}}%
    \special{pa 2610 1659}%
    \special{pa 2860 1225}%
    \special{fp}%
    \special{pa 2506 1579}%
    \special{pa 2860 1225}%
    \special{fp}%
    \special{pa 2426 1476}%
    \special{pa 2860 1225}%
    \special{fp}%
    \graphtemp=.5ex\advance\graphtemp by 1.001in
    \rlap{\kern 2.230in\lower\graphtemp\hbox to 0pt{\hss $\bullet$\hss}}%
    \graphtemp=.5ex\advance\graphtemp by 0.871in
    \rlap{\kern 2.213in\lower\graphtemp\hbox to 0pt{\hss $\bullet$\hss}}%
    \graphtemp=.5ex\advance\graphtemp by 0.742in
    \rlap{\kern 2.230in\lower\graphtemp\hbox to 0pt{\hss $\bullet$\hss}}%
    \special{pa 2230 1001}%
    \special{pa 2713 871}%
    \special{fp}%
    \special{pa 2213 871}%
    \special{pa 2713 871}%
    \special{fp}%
    \special{pa 2230 742}%
    \special{pa 2713 871}%
    \special{fp}%
    \graphtemp=.5ex\advance\graphtemp by 0.267in
    \rlap{\kern 2.426in\lower\graphtemp\hbox to 0pt{\hss $\bullet$\hss}}%
    \graphtemp=.5ex\advance\graphtemp by 0.163in
    \rlap{\kern 2.506in\lower\graphtemp\hbox to 0pt{\hss $\bullet$\hss}}%
    \graphtemp=.5ex\advance\graphtemp by 0.083in
    \rlap{\kern 2.610in\lower\graphtemp\hbox to 0pt{\hss $\bullet$\hss}}%
    \special{pa 2426 267}%
    \special{pa 2860 517}%
    \special{fp}%
    \special{pa 2506 163}%
    \special{pa 2860 517}%
    \special{fp}%
    \special{pa 2610 83}%
    \special{pa 2860 517}%
    \special{fp}%
    \graphtemp=.5ex\advance\graphtemp by 2.207in
    \rlap{\kern 3.381in\lower\graphtemp\hbox to 0pt{\hss (b)\hss}}%
    \hbox{\vrule depth2.207in width0pt height 0pt}%
    \kern 4.800in
  }%
}%
}

  \vspace{1pt}

  \caption{(a) Erd\H{o}s and Ne\v{s}et\v{r}il's construction for $\Delta=4$.  This graph requires 20 colors. (b) The largest possible neighborhood of an edge has size 24.}

\end{figure}

Let $v$ be an arbitrary vertex of a graph $G$.  Let $\mbox{dist}_v(v_1)$ denote the distance from vertex $v_1$ to $v$.  Let {\it distance class} $i$ be the set of vertices at distance $i$ from vertex $v$.  The distance class of an edge is the minimum of the distance classes of its vertices.  We say that an edge order is {\it compatible} with vertex $v$ if $e_1$ precedes $e_2$ in the order only when $\mbox{dist}_v(e_1)\geq \mbox{dist}_v(e_2)$.  Intuitively, we color all the edges in distance class $i+1$ (farther from $v$) before we color any edge in distance class $i$ (nearer to $v$).  Similarly, if we specify a cycle $C$ in the graph, we can define distance class $i$ to be the set of vertices at distance $i$ from cycle $C$.  We say an edge order is compatible with $C$ if $e_1$ precedes $e_2$ in the order only when $dist_C(e_1)\geq dist_C(e_2)$.
%\vspace{.4in}

\begin{lemma}If $G$ is a graph with maximum degree 4, then $G$ has a strong edge-coloring that uses 21 colors except that it leaves uncolored those edges incident to a single vertex.  If $C$ is a cycle in $G$, then $G$ has a strong edge-coloring that uses 21 colors except that it leaves uncolored the edges of $C$.
\label{simple}
\end{lemma}
\textbf{Proof.} We first consider the case of leaving uncolored only the edges incident to a single vertex.  Let $v$ be a vertex of $G$.  
Greedily color the edges in an order that is compatible with vertex $v$.  Suppose we are coloring edge $e$, not incident to $v$.  Let $u$ be a vertex adjacent to an endpoint of $e$ that is on a shortest path from $e$ to $v$.  Then none of the four edges incident to $u$ has been colored, since each edge incident to $u$ is in a lower distance class than $e$.  Thus, $|N(e)|\leq 24 - 4 = 20$.  

To prove the case of leaving uncolored only the edges of $C$, we color the edges in an order compatible with $C$.  The argument above holds for every edge not incident to $C$.  If $e$ is incident to $C$ and $|C|\geq 4$, then at least four edges in the neighborhood of $e$ are edges of $C$; so again $|N(e)|\leq 24 - 4 =20$.  If $e$ is incident to $C$ and $|C|=3$, then by counting we see that the neighborhood of $e$ has size at most 23.  The three uncolored edges of $C$ imply that $|N(e)|\leq 23 - 3= 20$.
$\Box$
\smallskip

Lemma \ref{simple} shows that if a graph has maximum degree 4 we can color nearly all edges using at most 21 colors.  In the rest of this paper, we show that we can always finish the edge-coloring using at most one additional color.
Theorem \ref{main1} is the main result of this paper.  We give the proof of the general case (4-regular and girth at least six) now, and defer the other cases (when the graph is not 4-regular or has small girth) to lemmas \ref{low degree}-\ref{5-cycle} in the remainder of the paper.

\begin{theorem}
Any graph with maximum degree 4 has a strong edge-coloring with at most 22 colors.
\label{main1}
\end{theorem}

\begin{lemma}
\label{general case}
Any 4-regular graph with girth at least six has a strong edge-coloring with at most 22 colors.
\end{lemma}
\textbf{Proof.}
By Lemma \ref{simple}, we choose an arbitrary vertex $v$, and greedily color all edges not incident to $v$, using at most 21 colors.  Now we recolor edges $e_1, e_2, e_3,$ and $e_4$ (as shown in figure 2) using color 22.  This allows us to greedily extend the coloring to the four edges incident to $v$.  Edges $e_1, e_2, e_3,$ and $e_4$ can receive the same color since the girth of $G$ is at least 6. 
$\Box$
\smallskip

\begin{figure}[hbtp]
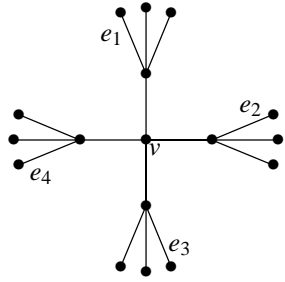

\gpic{
\expandafter\ifx\csname graph\endcsname\relax \csname newbox\endcsname\graph\fi
\expandafter\ifx\csname graphtemp\endcsname\relax \csname newdimen\endcsname\graphtemp\fi
\setbox\graph=\vtop{\vskip 0pt\hbox{%
    \special{pn 8}%
    \special{pa 750 404}%
    \special{pa 750 750}%
    \special{pa 1096 750}%
    \special{pa 750 750}%
    \special{pa 750 1096}%
    \special{pa 750 750}%
    \special{pa 404 750}%
    \special{fp}%
    \graphtemp=.5ex\advance\graphtemp by 0.750in
    \rlap{\kern 0.750in\lower\graphtemp\hbox to 0pt{\hss $\bullet$\hss}}%
    \graphtemp=.5ex\advance\graphtemp by 0.404in
    \rlap{\kern 0.750in\lower\graphtemp\hbox to 0pt{\hss $\bullet$\hss}}%
    \graphtemp=.5ex\advance\graphtemp by 0.750in
    \rlap{\kern 1.096in\lower\graphtemp\hbox to 0pt{\hss $\bullet$\hss}}%
    \graphtemp=.5ex\advance\graphtemp by 1.096in
    \rlap{\kern 0.750in\lower\graphtemp\hbox to 0pt{\hss $\bullet$\hss}}%
    \graphtemp=.5ex\advance\graphtemp by 0.750in
    \rlap{\kern 0.404in\lower\graphtemp\hbox to 0pt{\hss $\bullet$\hss}}%
    \graphtemp=.5ex\advance\graphtemp by 0.084in
    \rlap{\kern 0.618in\lower\graphtemp\hbox to 0pt{\hss $\bullet$\hss}}%
    \graphtemp=.5ex\advance\graphtemp by 0.058in
    \rlap{\kern 0.750in\lower\graphtemp\hbox to 0pt{\hss $\bullet$\hss}}%
    \graphtemp=.5ex\advance\graphtemp by 0.084in
    \rlap{\kern 0.882in\lower\graphtemp\hbox to 0pt{\hss $\bullet$\hss}}%
    \special{pa 618 84}%
    \special{pa 750 404}%
    \special{fp}%
    \special{pa 750 58}%
    \special{pa 750 404}%
    \special{fp}%
    \special{pa 882 84}%
    \special{pa 750 404}%
    \special{fp}%
    \graphtemp=.5ex\advance\graphtemp by 0.201in
    \rlap{\kern 0.572in\lower\graphtemp\hbox to 0pt{\hss $e_1$\hss}}%
    \graphtemp=.5ex\advance\graphtemp by 0.618in
    \rlap{\kern 1.416in\lower\graphtemp\hbox to 0pt{\hss $\bullet$\hss}}%
    \graphtemp=.5ex\advance\graphtemp by 0.750in
    \rlap{\kern 1.442in\lower\graphtemp\hbox to 0pt{\hss $\bullet$\hss}}%
    \graphtemp=.5ex\advance\graphtemp by 0.882in
    \rlap{\kern 1.416in\lower\graphtemp\hbox to 0pt{\hss $\bullet$\hss}}%
    \special{pa 1416 618}%
    \special{pa 1096 750}%
    \special{fp}%
    \special{pa 1442 750}%
    \special{pa 1096 750}%
    \special{fp}%
    \special{pa 1416 882}%
    \special{pa 1096 750}%
    \special{fp}%
    \graphtemp=.5ex\advance\graphtemp by 0.572in
    \rlap{\kern 1.299in\lower\graphtemp\hbox to 0pt{\hss $e_2$\hss}}%
    \graphtemp=.5ex\advance\graphtemp by 1.416in
    \rlap{\kern 0.882in\lower\graphtemp\hbox to 0pt{\hss $\bullet$\hss}}%
    \graphtemp=.5ex\advance\graphtemp by 1.442in
    \rlap{\kern 0.750in\lower\graphtemp\hbox to 0pt{\hss $\bullet$\hss}}%
    \graphtemp=.5ex\advance\graphtemp by 1.416in
    \rlap{\kern 0.618in\lower\graphtemp\hbox to 0pt{\hss $\bullet$\hss}}%
    \special{pa 882 1416}%
    \special{pa 750 1096}%
    \special{fp}%
    \special{pa 750 1442}%
    \special{pa 750 1096}%
    \special{fp}%
    \special{pa 618 1416}%
    \special{pa 750 1096}%
    \special{fp}%
    \graphtemp=.5ex\advance\graphtemp by 1.299in
    \rlap{\kern 0.928in\lower\graphtemp\hbox to 0pt{\hss $e_3$\hss}}%
    \graphtemp=.5ex\advance\graphtemp by 0.882in
    \rlap{\kern 0.084in\lower\graphtemp\hbox to 0pt{\hss $\bullet$\hss}}%
    \graphtemp=.5ex\advance\graphtemp by 0.750in
    \rlap{\kern 0.058in\lower\graphtemp\hbox to 0pt{\hss $\bullet$\hss}}%
    \graphtemp=.5ex\advance\graphtemp by 0.618in
    \rlap{\kern 0.084in\lower\graphtemp\hbox to 0pt{\hss $\bullet$\hss}}%
    \special{pa 84 882}%
    \special{pa 404 750}%
    \special{fp}%
    \special{pa 58 750}%
    \special{pa 404 750}%
    \special{fp}%
    \special{pa 84 618}%
    \special{pa 404 750}%
    \special{fp}%
    \graphtemp=.5ex\advance\graphtemp by 0.928in
    \rlap{\kern 0.201in\lower\graphtemp\hbox to 0pt{\hss $e_4$\hss}}%
    \graphtemp=.5ex\advance\graphtemp by 0.796in
    \rlap{\kern 0.796in\lower\graphtemp\hbox to 0pt{\hss $v$\hss}}%
    \hbox{\vrule depth1.500in width0pt height 0pt}%
    \kern 1.500in
  }%
}%
}
\caption{Vertex $v$ has degree 4 and the girth of the graph is at least 6.}
\end{figure}

Lemma \ref{general case} proves Theorem \ref{main1} for 4-regular graphs with girth at least 6.  To prove Theorem \ref{main1} for graphs that are not 4-regular and graphs with girth less than six, we use two ideas.  In Section \ref{small cases}, we consider graphs that are not 4-regular and graphs with girth at most 3.  In each case, we exploit local structure of the graph to give an edge order with $|N(e)|\leq 20$ for every edge in $G$.  In Sections \ref{girth four} and \ref{girth five}, we consider 4-regular graphs with girth 4 or 5.  We find pairs of edges that can receive the same color.  In this case, even though $|N(e)| > 21$, because not every edge in $N(e)$ receives a distinct color, we ensure that at most 22 colors are used.

\section{Graphs with $\delta<4$ or girth at most 3}
\label{small cases}

The three lemmas in this section are each proved using the same idea.  We color nearly all edges as in Lemma \ref{simple}.  We show that due to the presence of a low degree vertex, a loop, a double edge, or a 3-cycle, it is possible to order the remaining uncolored edges so that a greedy coloring uses at most 21 colors.

\begin{lemma}Any graph with maximum degree 4 that has a vertex with degree at most 3 has a strong edge-coloring that uses 21 colors.
\label{low degree}
\end{lemma}
\textbf{Proof.}  We assume $d(v)=3$ (if actually $d(v)<3$, this only makes it easier to complete the coloring).  Color the edges in an order that is compatible with vertex $v$.  Let $e_1,e_2,e_3$ be the edges incident to vertex $v$.  If the edges are ordered $e_1$, $e_2$, $e_3$, we have ${|N(e_1)|\leq 18}, {|N(e_2)|\leq 19}$ and ${|N(e_3)|\leq 20}$, so there are colors for $e_1$, $e_2$ and $e_3$. 
$\Box$

\begin{lemma}A 4-regular graph with a loop or a double edge has a strong edge-coloring that uses 21 colors.
\end{lemma}
\textbf{Proof.} If $e$ is a loop incident to vertex $v$, we can greedily color the edges in an order compatible with vertex $v$ (this is very similar to Lemma \ref{low degree}).  So we can assume that the graph has a double edge.

Let $v$ be one of the vertices incident to the double edge.  Color the edges in an order that is compatible with vertex $v$.  Let $e_3,e_4$ be the double edges and $e_1,e_2$ be the other edges incident to $v$. Then ${|N(e_1)|\leq 17}, {|N(e_2)|\leq 18}, {|N(e_3)|\leq 16}$ and ${|N(e_4)|\leq 17}$, so there are colors for $e_1$, $e_2$, $e_3$, and $e_4$. $\Box$

\begin{lemma}A 4-regular graph with girth 3 has a strong edge-coloring that uses 21 colors.
\end{lemma}
\textbf{Proof.} Let $C$ be a 3-cycle in the graph.  By Lemma \ref{simple} we greedily color all edges except the edges of $C$; this uses at most 21 colors.  An edge of a 3-cycle has a neighborhood of size at most 20, so each of the three uncolored edges satisfies $|N(e)|\leq 18$ and we can greedily finish the coloring.
$\Box$

\section{4-regular graphs with girth four}
\label{girth four}
Lemma \ref{simple} shows that we can color nearly all edges of the graph using 21 colors.  Here we consider 4-regular graphs of girth four.  We give an edge order such that the greedy coloring uses at most 22 colors; in some cases we precolor four edges prior to the greedy coloring.  We use $A(e)$ to denote the set of colors available on edge $e$.

%lemma 7
\begin{lemma}Any 4-regular graph with girth 4 has a strong edge-coloring that uses 22 colors.
\label{4-cycle}
\end{lemma}
\textbf{Proof.}
Let $C$ be a 4-cycle, with the 4 edges labeled $c_i$ ($1\leq i\leq 4$) in
clockwise order and the pair of edges not on the cycle and adjacent to $c_i$ and $c_{i-1}$ is labeled $a_i$ and $b_i$ (all subscripts are $\mbox{mod }4$).  
 We refer to the edges labeled by $a_i$ and $b_i$ as {\it incident} edges.  
By Lemma \ref{simple}, we greedily color all edges except the edges of $C$ and the 8 incident edges.  This uses at most 21 colors.  If two incident edges share an endpoint not on $C$, the two edges form an {\it adjacent pair}.  The only possibility of an adjacent pair is if $a_1$ or $b_1$ shares an endpoint with $a_3$ or $b_3$ (or similarly if $a_2$ or $b_2$ shares an endpoint with $a_4$ or $b_4$).  If the twelve uncolored edges contain at least two adjacent pairs, then we greedily color the incident edges.  The neighborhood of each $c_i$ has size at most 21, so $|A(c_i)|\geq 4$ for all $i$; thus we can finish by greedily coloring the four edges of $C$.  

\begin{figure}[hbtp]
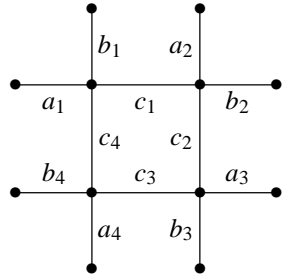

\gpic{
\expandafter\ifx\csname graph\endcsname\relax \csname newbox\endcsname\graph\fi
\expandafter\ifx\csname graphtemp\endcsname\relax \csname newdimen\endcsname\graphtemp\fi
\setbox\graph=\vtop{\vskip 0pt\hbox{%
    \graphtemp=.5ex\advance\graphtemp by 0.467in
    \rlap{\kern 1.033in\lower\graphtemp\hbox to 0pt{\hss $\bullet$\hss}}%
    \graphtemp=.5ex\advance\graphtemp by 1.033in
    \rlap{\kern 1.033in\lower\graphtemp\hbox to 0pt{\hss $\bullet$\hss}}%
    \graphtemp=.5ex\advance\graphtemp by 1.033in
    \rlap{\kern 0.467in\lower\graphtemp\hbox to 0pt{\hss $\bullet$\hss}}%
    \graphtemp=.5ex\advance\graphtemp by 0.467in
    \rlap{\kern 0.467in\lower\graphtemp\hbox to 0pt{\hss $\bullet$\hss}}%
    \special{pn 8}%
    \special{pa 1033 467}%
    \special{pa 1033 1033}%
    \special{pa 467 1033}%
    \special{pa 467 467}%
    \special{pa 1033 467}%
    \special{fp}%
    \special{pa 1433 467}%
    \special{pa 1033 467}%
    \special{pa 1033 67}%
    \special{fp}%
    \graphtemp=.5ex\advance\graphtemp by 0.467in
    \rlap{\kern 1.433in\lower\graphtemp\hbox to 0pt{\hss $\bullet$\hss}}%
    \graphtemp=.5ex\advance\graphtemp by 0.067in
    \rlap{\kern 1.033in\lower\graphtemp\hbox to 0pt{\hss $\bullet$\hss}}%
    \special{pa 1033 1433}%
    \special{pa 1033 1033}%
    \special{pa 1433 1033}%
    \special{fp}%
    \graphtemp=.5ex\advance\graphtemp by 1.433in
    \rlap{\kern 1.033in\lower\graphtemp\hbox to 0pt{\hss $\bullet$\hss}}%
    \graphtemp=.5ex\advance\graphtemp by 1.033in
    \rlap{\kern 1.433in\lower\graphtemp\hbox to 0pt{\hss $\bullet$\hss}}%
    \special{pa 67 1033}%
    \special{pa 467 1033}%
    \special{pa 467 1433}%
    \special{fp}%
    \graphtemp=.5ex\advance\graphtemp by 1.033in
    \rlap{\kern 0.067in\lower\graphtemp\hbox to 0pt{\hss $\bullet$\hss}}%
    \graphtemp=.5ex\advance\graphtemp by 1.433in
    \rlap{\kern 0.467in\lower\graphtemp\hbox to 0pt{\hss $\bullet$\hss}}%
    \special{pa 467 67}%
    \special{pa 467 467}%
    \special{pa 67 467}%
    \special{fp}%
    \graphtemp=.5ex\advance\graphtemp by 0.067in
    \rlap{\kern 0.467in\lower\graphtemp\hbox to 0pt{\hss $\bullet$\hss}}%
    \graphtemp=.5ex\advance\graphtemp by 0.467in
    \rlap{\kern 0.067in\lower\graphtemp\hbox to 0pt{\hss $\bullet$\hss}}%
    \graphtemp=.5ex\advance\graphtemp by 0.563in
    \rlap{\kern 0.750in\lower\graphtemp\hbox to 0pt{\hss $c_1$\hss}}%
    \graphtemp=.5ex\advance\graphtemp by 0.750in
    \rlap{\kern 0.937in\lower\graphtemp\hbox to 0pt{\hss $c_2$\hss}}%
    \graphtemp=.5ex\advance\graphtemp by 0.937in
    \rlap{\kern 0.750in\lower\graphtemp\hbox to 0pt{\hss $c_3$\hss}}%
    \graphtemp=.5ex\advance\graphtemp by 0.750in
    \rlap{\kern 0.563in\lower\graphtemp\hbox to 0pt{\hss $c_4$\hss}}%
    \graphtemp=.5ex\advance\graphtemp by 0.563in
    \rlap{\kern 0.270in\lower\graphtemp\hbox to 0pt{\hss $a_1$\hss}}%
    \graphtemp=.5ex\advance\graphtemp by 0.270in
    \rlap{\kern 0.937in\lower\graphtemp\hbox to 0pt{\hss $a_2$\hss}}%
    \graphtemp=.5ex\advance\graphtemp by 0.937in
    \rlap{\kern 1.230in\lower\graphtemp\hbox to 0pt{\hss $a_3$\hss}}%
    \graphtemp=.5ex\advance\graphtemp by 1.230in
    \rlap{\kern 0.563in\lower\graphtemp\hbox to 0pt{\hss $a_4$\hss}}%
    \graphtemp=.5ex\advance\graphtemp by 0.270in
    \rlap{\kern 0.563in\lower\graphtemp\hbox to 0pt{\hss $b_1$\hss}}%
    \graphtemp=.5ex\advance\graphtemp by 0.563in
    \rlap{\kern 1.230in\lower\graphtemp\hbox to 0pt{\hss $b_2$\hss}}%
    \graphtemp=.5ex\advance\graphtemp by 1.230in
    \rlap{\kern 0.937in\lower\graphtemp\hbox to 0pt{\hss $b_3$\hss}}%
    \graphtemp=.5ex\advance\graphtemp by 0.937in
    \rlap{\kern 0.270in\lower\graphtemp\hbox to 0pt{\hss $b_4$\hss}}%
    \hbox{\vrule depth1.500in width0pt height 0pt}%
    \kern 1.500in
  }%
}%
}
\caption{A 4-cycle in a 4-regular graph.}
\end{figure}

Suppose the uncolored edges contain exactly one adjacent pair.  For example, suppose edges $a_2$ and $a_4$ share an endpoint.  Call edges $a_1,b_1,a_3$, and $b_3$ a {\it pack}.  Consider the case when we can assign color 22 to two edges of the pack.  Now we greedily color all edges except the edges of $C$.  This uses at most 21 colors (Lemma \ref{simple}).  Each $c_i$ has a neighborhood with size at most 22.  Since color 22 is used twice in the neighborhood of each $c_i$, each $c_i$ satisfies $|A(c_i)|\geq 4$.  So we can greedily finish the coloring.
Instead consider the case when no pair of edges in the pack can receive the same color.  This implies the existence of edges between each pair of nonadjacent edges of the pack.  Call these four additional edges {\it diagonal} edges.  Observe (by counting) that the neighborhood of a diagonal edge has size at most 21.  So we can color the diagonal edges last in the greedy coloring.  Thus we greedily color all edges except the four edges of $C$ and the four diagonal edges (this uses at most 21 colors).  Now we color the four edges of $C$ (the four uncolored diagonal edges ensure there are enough colors available to color the edges of $C$).  Lastly, we color the four diagonal edges.

Finally, suppose that the uncolored edges contain no adjacent pairs.  In this case we will greedily color almost all edges of the graph (Lemma \ref{simple}), but must do additional work beforehand to ensure that after greedily coloring most of the edges each $c_i$ will satisfy $|A(c_i)|\geq 4$.  
As above, call edges $a_1,b_1,a_3$, and $b_3$ a pack.  Similarly, call edges $a_2,b_2,a_4$, and $b_4$ a pack.  

Consider the case when we can assign color 21 to two edges of one pack and assign color 22 to two edges of the other pack.  We greedily color all edges but the four edges of $C$.  
  Lemma \ref{simple} showed that a similar greedy coloring used at most 21 colors; however in Lemma \ref{simple} none of the edges were precolored.  We adapt that argument to show that even in the presence of these four precolored edges a greedy coloring uses at most 22 colors.  Lemma \ref{simple} argued there were at least four uncolored edges in the neighborhood of the edge being colored, so $|N(e)|\leq 20$.  The same argument applies in this case except that possibly one of the edges that was uncolored in Lemma \ref{simple} is now colored.  Hence $|N(e)|\leq 21$ (this follows from the fact that the four uncolored edges in Lemma \ref{simple} were incident to the same vertex and in the present situation at most one precolored edge is incident to each vertex). Hence, the greedy coloring uses at most 22 colors.
The neighborhood of each $c_i$ has size at most 23.  Since colors 21 and 22 are each repeated in the neighborhood of each $c_i$, we see that each $c_i$ satisfies $|A(c_i)|\geq 4$.  So we can greedily finish the coloring.

Instead, consider the case when we can not assign color 21 to two edges of one pack and assign color 22 to two edges of the other pack.  If no two edges in a pack can receive the same color, this implies the existence of edges between each pair of nonadjacent edges of the pack.  These four diagonal edges each have a neighborhood with size at most 21.  As we did above, we greedily color all edges except the four edges of $C$ and the four diagonal edges.  Now we color the four edges of $C$, and lastly, we color the four diagonal edges.
$\Box$

\section{4-regular graphs with girth five}
\label{girth five}

Here we consider 4-regular graphs with girth five.  As in the case of girth four, we color nearly all the edges by Lemma \ref{simple}.
Intuitively, if there are enough different colors available to be used on the remaining uncolored edges, we should be able to complete this coloring by giving each uncolored edge its own color.  However, if there are fewer different colors available than the number of uncolored edges, this approach is doomed to fail.  Hall's Theorem formalizes this intuition.  In the language of Hall's Theorem, we have $m$ uncolored edges, and the set $A_i$ denotes the colors available to use on edge $i$.  For a proof of Hall's Theorem, we refer the reader to {\it Introduction to Graph Theory} \cite{We3}.

\begin{theorem}[Hall's Theorem]
\label{Hall's Theorem}
There exists a system of distinct representatives for a family of sets $A_1, A_2 , \ldots, A_m$ if and only if the union of any $j$ of these sets contains at least $j$ elements for all $j$ from $1$ to $m$.
\end{theorem}

We define a {\it partial coloring} to be a strong edge-coloring except that some edges may be uncolored.
Suppose that we have a partial coloring, with only the edge set $T$ left uncolored.  Let $A(e)$ be the set of colors available to color edge $e$.  Then Hall's Theorem guarantees that if we are unable to complete the coloring by giving each edge its own color, there exists a set $S \subseteq T$ with $|S| > |\cup_{e\in S}A(e)|$.  Define the \textit{discrepancy}, $disc(S) = |S| - |\cup_{e\in S}A(e)|$.

Our idea is to color the set of edges with maximum discrepancy, then argue that this coloring can be extended to the remaining uncolored edges.

%lemma 9
\begin{lemma}
\label{Hall's Theorem Extension}
Let $T$ be the set of uncolored edges in a partially colored graph.  Let $S$ be a subset of $T$ with maximum discrepancy.  Then a valid coloring for $S$ can be extended to a valid coloring for all of $T$.
\end{lemma}
\textbf{Proof.}  
%Assume the claim is false.  Since the coloring of $S$ cannot be extended to $T\setminus S$, some set of edges $S'\subseteq (T\setminus S)$ has positive discrepancy (after coloring $S$).  We show that $disc(S\cup S') > disc(S)$.  Let $k = disc(S)$ and $l = disc(S')$ (after coloring $S$).  This gives the equations $|S| = k + |\cup_{e\in S}A(e)|$ and $|S'| = l + |\cup_{e\in S'}A(e)| - |(\cup_{e\in S}A(e)) \cap (\cup_{e\in S'}A(e))|$.  Adding these equations, we get $|S\cup S'| = |S| + |S'| = k + l + |\cup_{e\in S}A(e)| + |\cup_{e\in S'}A(e)| - |(\cup_{e\in S}A(e)) \cap (\cup_{e\in S'}A(e))| = k + l + |\cup_{e\in S\cup S'}A(e)|$.  Thus $disc(S\cup S') = k + l > k = disc(S)$.
%
Assume the claim is false.  Since the coloring of $S$ cannot be extended to $T\setminus S$, some set of edges $S' \subseteq (T\setminus S)$ has positive discrepancy (after coloring $S$).  We show that $\mbox{disc}(S\cup S') > \mbox{disc}(S)$.  Let $R$ be the set of colors available to use on at least one edge of $(S \cup S')$.  Let $R_1$ be the set of colors available to use on at least one edge of $S$.  Let $R_2$ be the set of colors available to use on at least one edge of $S'$ after the edges of $S$ have been colored.  Let $k = disc(S)$.  Then
$|S| = k + |R_1|$ and $|S'| \geq 1 + |R_2|$.  Since $S$ and $S'$ are disjoint, we get
$$
       |S \cup S'| = |S| + |S'| \geq k + 1 + |R_1| + |R_2| > k + |R|.
$$

The latter inequality holds since a color which is in $R\setminus R_1$ must be in $R_2$ and therefore we have $|R|=|R_1\cup R_2|\leq |R_1| + |R_2|$.  Hence
$$
       \mbox{disc}(S\cup S') = |S\cup S'| - |R| > k = \mbox{disc}(S)
$$

This contradicts the maximality of disc($S$).  Hence, any valid coloring of $S$ can be extended to a valid coloring of $T$. 
$\Box$
\smallskip

%lemma 10
\begin{lemma}If $G$ is a 4-regular graph with girth 5, then $G$ has a strong edge-coloring that uses 22 colors.
\label{5-cycle}
\end{lemma}
\textbf{Proof.} 
Let $C$ be a 5-cycle, with the 5 edges labeled $c_i$ ($1\leq i\leq 5$) in
clockwise order and the pair of edges not on the cycle and adjacent to $c_i$ and $c_{i-1}$ is labeled $a_i$ and $b_i$ (all subscripts are $\mbox{mod }5$).  We refer to the edges labeled by $a_i$ and $b_i$ as {\it incident} edges.  Edge $a_1$ is at least distance 2 from at least one of edges $a_3$ and $b_3$; for if $a_1$ has edge $e_1$ to $a_3$ and edge $e_2$ to $b_3$ then we have the 4-cycle $e_1,e_2,b_3,a_3$.  Thus (by possibly renaming $a_3$ and $b_3$) we can assume there is no edge between edges $a_1$ and $b_3$.

\begin{figure}[hbtp]
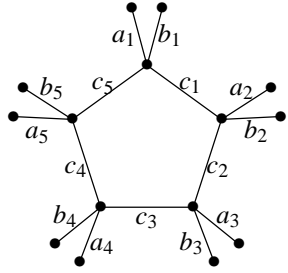

\gpic{
\expandafter\ifx\csname graph\endcsname\relax \csname newbox\endcsname\graph\fi
\expandafter\ifx\csname graphtemp\endcsname\relax \csname newdimen\endcsname\graphtemp\fi
\setbox\graph=\vtop{\vskip 0pt\hbox{%
    \graphtemp=.5ex\advance\graphtemp by 0.349in
    \rlap{\kern 0.750in\lower\graphtemp\hbox to 0pt{\hss $\bullet$\hss}}%
    \graphtemp=.5ex\advance\graphtemp by 0.633in
    \rlap{\kern 1.141in\lower\graphtemp\hbox to 0pt{\hss $\bullet$\hss}}%
    \graphtemp=.5ex\advance\graphtemp by 1.093in
    \rlap{\kern 0.992in\lower\graphtemp\hbox to 0pt{\hss $\bullet$\hss}}%
    \graphtemp=.5ex\advance\graphtemp by 1.093in
    \rlap{\kern 0.508in\lower\graphtemp\hbox to 0pt{\hss $\bullet$\hss}}%
    \graphtemp=.5ex\advance\graphtemp by 0.633in
    \rlap{\kern 0.359in\lower\graphtemp\hbox to 0pt{\hss $\bullet$\hss}}%
    \special{pn 8}%
    \special{pa 750 349}%
    \special{pa 1141 633}%
    \special{pa 992 1093}%
    \special{pa 508 1093}%
    \special{pa 359 633}%
    \special{pa 750 349}%
    \special{fp}%
    \special{pa 830 51}%
    \special{pa 750 349}%
    \special{pa 670 51}%
    \special{fp}%
    \graphtemp=.5ex\advance\graphtemp by 0.051in
    \rlap{\kern 0.830in\lower\graphtemp\hbox to 0pt{\hss $\bullet$\hss}}%
    \graphtemp=.5ex\advance\graphtemp by 0.051in
    \rlap{\kern 0.670in\lower\graphtemp\hbox to 0pt{\hss $\bullet$\hss}}%
    \special{pa 1449 617}%
    \special{pa 1141 633}%
    \special{pa 1399 465}%
    \special{fp}%
    \graphtemp=.5ex\advance\graphtemp by 0.617in
    \rlap{\kern 1.449in\lower\graphtemp\hbox to 0pt{\hss $\bullet$\hss}}%
    \graphtemp=.5ex\advance\graphtemp by 0.465in
    \rlap{\kern 1.399in\lower\graphtemp\hbox to 0pt{\hss $\bullet$\hss}}%
    \special{pa 1102 1380}%
    \special{pa 992 1093}%
    \special{pa 1231 1286}%
    \special{fp}%
    \graphtemp=.5ex\advance\graphtemp by 1.380in
    \rlap{\kern 1.102in\lower\graphtemp\hbox to 0pt{\hss $\bullet$\hss}}%
    \graphtemp=.5ex\advance\graphtemp by 1.286in
    \rlap{\kern 1.231in\lower\graphtemp\hbox to 0pt{\hss $\bullet$\hss}}%
    \special{pa 269 1286}%
    \special{pa 508 1093}%
    \special{pa 398 1380}%
    \special{fp}%
    \graphtemp=.5ex\advance\graphtemp by 1.286in
    \rlap{\kern 0.269in\lower\graphtemp\hbox to 0pt{\hss $\bullet$\hss}}%
    \graphtemp=.5ex\advance\graphtemp by 1.380in
    \rlap{\kern 0.398in\lower\graphtemp\hbox to 0pt{\hss $\bullet$\hss}}%
    \special{pa 101 465}%
    \special{pa 359 633}%
    \special{pa 51 617}%
    \special{fp}%
    \graphtemp=.5ex\advance\graphtemp by 0.465in
    \rlap{\kern 0.101in\lower\graphtemp\hbox to 0pt{\hss $\bullet$\hss}}%
    \graphtemp=.5ex\advance\graphtemp by 0.617in
    \rlap{\kern 0.051in\lower\graphtemp\hbox to 0pt{\hss $\bullet$\hss}}%
    \graphtemp=.5ex\advance\graphtemp by 0.197in
    \rlap{\kern 0.630in\lower\graphtemp\hbox to 0pt{\hss $a_1$\hss}}%
    \graphtemp=.5ex\advance\graphtemp by 0.197in
    \rlap{\kern 0.870in\lower\graphtemp\hbox to 0pt{\hss $b_1$\hss}}%
    \graphtemp=.5ex\advance\graphtemp by 0.444in
    \rlap{\kern 0.979in\lower\graphtemp\hbox to 0pt{\hss $c_1$\hss}}%
    \graphtemp=.5ex\advance\graphtemp by 0.472in
    \rlap{\kern 1.248in\lower\graphtemp\hbox to 0pt{\hss $a_2$\hss}}%
    \graphtemp=.5ex\advance\graphtemp by 0.700in
    \rlap{\kern 1.322in\lower\graphtemp\hbox to 0pt{\hss $b_2$\hss}}%
    \graphtemp=.5ex\advance\graphtemp by 0.881in
    \rlap{\kern 1.121in\lower\graphtemp\hbox to 0pt{\hss $c_2$\hss}}%
    \graphtemp=.5ex\advance\graphtemp by 1.145in
    \rlap{\kern 1.178in\lower\graphtemp\hbox to 0pt{\hss $a_3$\hss}}%
    \graphtemp=.5ex\advance\graphtemp by 1.286in
    \rlap{\kern 0.984in\lower\graphtemp\hbox to 0pt{\hss $b_3$\hss}}%
    \graphtemp=.5ex\advance\graphtemp by 1.150in
    \rlap{\kern 0.750in\lower\graphtemp\hbox to 0pt{\hss $c_3$\hss}}%
    \graphtemp=.5ex\advance\graphtemp by 1.286in
    \rlap{\kern 0.516in\lower\graphtemp\hbox to 0pt{\hss $a_4$\hss}}%
    \graphtemp=.5ex\advance\graphtemp by 1.145in
    \rlap{\kern 0.322in\lower\graphtemp\hbox to 0pt{\hss $b_4$\hss}}%
    \graphtemp=.5ex\advance\graphtemp by 0.881in
    \rlap{\kern 0.379in\lower\graphtemp\hbox to 0pt{\hss $c_4$\hss}}%
    \graphtemp=.5ex\advance\graphtemp by 0.700in
    \rlap{\kern 0.178in\lower\graphtemp\hbox to 0pt{\hss $a_5$\hss}}%
    \graphtemp=.5ex\advance\graphtemp by 0.472in
    \rlap{\kern 0.252in\lower\graphtemp\hbox to 0pt{\hss $b_5$\hss}}%
    \graphtemp=.5ex\advance\graphtemp by 0.444in
    \rlap{\kern 0.521in\lower\graphtemp\hbox to 0pt{\hss $c_5$\hss}}%
    \hbox{\vrule depth1.432in width0pt height 0pt}%
    \kern 1.500in
  }%
}%
}
\caption{A 5-cycle in a 4-regular graph.}
\end{figure}

By repeating the same argument, we can assume there is no edge between the two edges of each of the following pairs: $(a_1,b_3)$, $(a_3,b_5)$, $(a_5,b_2)$, and $(a_2,b_4)$.  Assign color 21 to edges $b_1$ and $c_3$ and assign color 22 to edges $a_5$ and $b_2$.  Greedily color all edges except the edges of $C$ and the incident edges.  This uses at most 22 colors. 
%(this is the same argument as in the last case of Lemma \ref{4-cycle}). 

There are 11 uncolored edges; if we can not assign a distinct color to each uncolored edge, then Hall's Theorem guarantees there exists a subset of the uncolored edges with positive discrepancy.  Let $S$ be a subset of the uncolored edges with maximum discrepancy.  By counting the uncolored edges in the neighborhood of each edge, we observe that if $e$ is an edge of $C$, then $|A(e)|\geq 8$ and if $e$ is an incident edge then $|A(e)|\geq 5$.  We can assume that $S$ contains some edge of $C$, since otherwise we can greedily color $S$ (Lemma \ref{simple}), then extend the coloring to the remaining uncolored edges (Lemma \ref{Hall's Theorem Extension}).  Since $disc(S)>0$ and $|A(e)|\geq 8$ for each edge of $C$, we have $|S|$ is 9, 10, or 11.

Suppose $|S|$ is 9 or 10.  
Then since $S$ is missing at most two uncolored edges, $S$ contains at least one of the pair $(a_1,b_3)$, the pair $(a_2,b_4)$, and the pair $(a_3,b_5)$.  
%By symmetry we can map the pair $(a_2,b_4)$ to the pair $(a_3,b_5)$, so we can assume $S$ contains either the pair $(a_1,b_3)$ or the pair $(a_3,b_5)$.  
Since each edge in the pair satisfies $|A(e)|\geq 5$ and $|\cup_{e\in S}A(e)|\leq 9$, some color is available for use on both edges of the pair.  
Assign the same color to both edges.  
Since the neighborhood of each uncolored incident edge, $e$, contains at least three uncolored edges of $C$, we have $|N(e)|\leq 24 - 3=21$; so we can greedily color the remaining uncolored incident edges.  
%Now color the edges of the 5-cycle in the order $c_2, c_4, c_5, c_1$.
Now if $S$ contains the pair $(a_1,b_3)$ or the pair $(a_3,b_5)$ then color the edges of the 5-cycle in the order $c_2,c_4,c_5,c_1$; if $S$ contains the pair $(a_2,b_4)$ then color the edges of the 5-cycle in the order $c_2,c_4,c_1,c_5$.

Suppose $|S|$ is 11 and that no color is available on both edges of any of the pairs $(a_1,b_3)$, $(a_2,b_4)$, and $(a_3,b_5)$ (otherwise the above argument holds).  Assign the same color to $c_1$ and $a_4$; call it color $x$.  Note that if $|A(c_1)|\geq 8$, $|A(a_4)|\geq 5$, and $|A(c_1)\cup A(a_4)|\leq |\cup_{e\in S}A(e)|\leq 10$, then $|A(c_1)\cap A(a_4)|\neq 0$.  Before color $x$ was assigned to $c_1$ and $a_4$, it had been available on exactly one edge of each of the three pairs.  Greedily color those three edges (none of the colors used on these three edges is color $x$).  Now the three remaining uncolored incident edges each satisfy $|A(e)|\geq 3$, so we can greedily color them.  Greedily color the three remaining edges in the order $c_2, c_4, c_5$.
$\Box$

\section{Conclusion}
We note that it is straightforward to convert this proof to an algorithm that runs in linear time.  We assume a data structure that stores all the relevant information about each vertex.  Using breadth-first search, we can calculate the distance classes, as well as implement each lemma in linear time.

A natural question is whether it is possible to extend the ideas of this paper to larger $\Delta$.  The best bound we could hope for from the techniques of this paper is $2{\Delta}^2 - 3{\Delta} + 2$.  It is straightforward to prove an analog of Lemma \ref{simple} that gives a strong edge-coloring of $G$ that uses $2{\Delta}^2-3{\Delta}+1$ colors except that it leaves uncolored those edges incident to a single vertex (however, the author was unable to prove an analog to the ``uncolored cycle'' portion of Lemma \ref{simple}).
If $G$ contains a loop, a double edge, or a vertex of degree less than ${\Delta}$, then by the analog of Lemma \ref{simple} $G$ has a strong edge-coloring that uses at most $2{\Delta}^2-3{\Delta}+1$ colors.  Using the ideas of Lemma \ref{general case}, we see that if $G$ is ${\Delta}$-regular and has girth at least 6, then $G$ has a strong edge-coloring that uses $2{\Delta}^2-3{\Delta}+2$ colors.  Thus, to complete a proof for graphs with larger $\Delta$, one must consider the case of regular graphs with girth 3, 4, or 5.

\section{Acknowledgements}

This paper draws heavily on ideas from a paper by Lars Andersen \cite{An}, in which he considers the case $\Delta(G)=3$.  The present author would like to thank him for that paper.  Without it, this one would not have been written.  The exposition of this paper has been greatly improved by critique from David Bunde and Erin Chambers.
\bibliography{strong-coloring}
\end{document}